\documentclass[12pt, oneside]{amsart}
\usepackage{amscd}
\usepackage{amssymb}
\usepackage{amsfonts}
\usepackage{remark}
\usepackage{enumerate}
\usepackage{xypic}

\newtheorem{theorem}{Theorem}[section]

\newtheorem{lemma}[theorem]{Lemma}
\newtheorem{proposition}[theorem]{Proposition}

\theoremstyle{definition}

\newtheorem{remark}[theorem]{Remark}
\newtheorem{emp}[theorem]{}

\theoremstyle{remark}

\newcommand{\N}{\mathcal N}
  
\newcommand{\Disc}{\mbox{\rm Disc}\kern 1pt}
\newcommand{\disc}{\mbox{\rm disc}\kern 1pt}
\newcommand{\Jac}{\mbox{\rm Jac}\kern 1pt}
\newcommand{\Res}{\mbox{\rm Res}\kern 1pt}
\newcommand{\Spec}{\mbox{\rm Spec}\kern 1pt}
\newcommand{\Br}{\mbox{\rm Br}\kern 1pt}
\newcommand{\Sch}{\mbox{\rm Sch}\kern 1pt}
\newcommand{\Sets}{\mbox{\rm Sets}\kern 1pt}

\newcommand{\Z}{\mathbb Z}

\newcommand{\Div}{\operatorname{Div}}
\newcommand{\Pic}{\operatorname{Pic}}

\newcommand{\Ker}{\operatorname{Ker}\kern 1pt}
\newcommand{\Coker}{\operatorname{Coker}\kern 1pt}
\newcommand{\Num}{\mbox{\rm Num}\kern 1pt}
\newcommand{\Image}{\mbox{\rm Im}\kern 1pt}

\newcommand{\Gal}{\mbox{\rm Gal}\kern 1pt}

\newcommand{\lcm}{\operatorname{lcm}}
\newcommand{\Lie}{\operatorname{Lie}\kern 1pt}
\newcommand{\LLie}{{\mathcal L}\mbox{\rm ie}\kern 1pt}
\newcommand{\chara}{\mbox{\rm char}\kern 1pt}

\DeclareFontEncoding{OT2}{}{} % to enable usage of cyrillic fonts
  \newcommand{\textcyr}[1]{%
    {\fontencoding{OT2}\fontfamily{wncyr}\fontseries{m}\fontshape{n}%
     \selectfont #1}}

\newcommand{\Sha}{{\mbox{\textcyr{Sh}}}}

\DeclareFontFamily{U}{wncy}{}
\DeclareFontShape{U}{wncy}{m}{n}{<->wncyr10}{}
\DeclareSymbolFont{mcy}{U}{wncy}{m}{n}
\DeclareMathSymbol{\Sh}{\mathord}{mcy}{"58}
\numberwithin{equation}{section}
\numberwithin{figure}{section}
\numberwithin{table}{section}
\if false
\makeatletter
    \let\c@equation\c@thm
    \let\c@figure\c@thm
   \let\c@table\c@thm
\makeatother
\fi
\begin{document}
\title[Corrigendum]{Corrigendum \\
to \\
 N\'eron models, Lie algebras, and reduction of curves of genus one \cite{LLR1}\\
 and \\
 The Brauer group of a surface \cite{LLR2}}
 
\author{Qing Liu}
\address{Institut de Math\'ematiques de Bordeaux\\
CNRS UMR 5251, Universit\'e de Bordeaux\\
33405 Talence cedex, France} 
\email{qing.liu@math.u-bordeaux.fr} 

\author{Dino Lorenzini}
\address{Department of Mathematics\\
University of Georgia\\
Athens, GA 30602, USA}
\email{lorenzin@uga.edu} 

\author{Michel Raynaud${}^{\dagger}$}\thanks{Our coauthor, mentor, and friend Michel Raynaud fell ill 
soon after we started writing this corrigendum. We are profoundly sad by his passing on March 10, 2018. All mistakes in this corrigendum are ours only (Qing Liu and Dino Lorenzini).} 
\address{Laboratoire de Math\'ematiques\\
Universit\'e Paris-Sud\\ 
91405 Orsay Cedex, FRANCE} 

\maketitle

\begin{section}{Introduction}

Let $k$ be a finite field of characteristic $p$. Let $V/k$ be a smooth projective 
geometrically connected curve with function field $K$. Let 
$X/k$ be a proper smooth and geometrically connected surface
endowed with a proper flat map $f : X\to V$ such that the 
generic fiber $X_K/K$ is smooth and geometrically connected of genus $g \geq 1$. 
Let $A_K/K$ denote the Jacobian of $X_K/K$. 
 
The proof of Theorem 4.3 in \cite{LLR1}, which we state in corrected form below, is based in part on a result of Gordon \cite{Gor}. 
Thomas Geisser noted in \cite{Gei} that the formula provided in Theorem 4.3 in \cite{LLR1} needs to be corrected, due to the fact 
that Lemma 4.2 in \cite{Gor} is 
missing a hypothesis. He provides a corrected formula in \cite{Gei}, Theorem 1.1, and his method applies also to the number field case
(up to a power of $2$ if not totally imaginary). 
Several of the
intermediate results in \cite{Gor}
are only valid 
under the assumption  
that $\Pic^0(X_K) = A_K(K)$. We revisit
the paper \cite{Gor} in this corrigendum to remove this hypothesis in all arguments.
In doing so, we also avoid using Lemma 4.3 in \cite{Gor}, whose proof is incorrect, and whose statement might be wrong in general.
\end{section}

\begin{section}{Corrected Statements}

We start by recalling the notation needed to state our main theorem. 
Let $\Sha(A_K)$ denote the Shafarevich-Tate
group of the abelian variety $A_K/K$. Let $\Br(X)$ denote the Brauer group of $X$.
It is well-known that if either $\Sha(A_K)$ or 
$\Br(X)$ is finite, then so is the other (see \cite{Tat}, section 3,  or \cite{GroB}, section 4).

 The {\em index} $\delta:=\delta(X_K)$ of a curve over a field $K$
is the least positive degree of a divisor on $X_K$.  
The {\em period} $\delta':=\delta'(X_K)$ of  $X_K$ is   
the order of the cokernel of the degree map $\Pic_{X_K/K}(K) \to {\mathbb Z}$.
When $v \in V $ is a closed point, we denote by $K_v$ the completion of $K$ at $v$, and let
$ \delta_v:= \delta(X_{K_v})$, and $\delta_v':= \delta'(X_{K_v})$. 

Recall that we have an exact sequence 
$$ 0 \longrightarrow \Pic^0(X_K) \longrightarrow A_K(K) \longrightarrow {\rm Br}(K).$$
Since the Brauer group ${\rm Br}(K)$ is a torsion group, and since
$A_K(K)$ is a finitely generated abelian group, 
the quotient $ A_K(K)/\Pic^0(X_K)$ is finite, and 
$\Pic^0(X_K)$ and $ A_K(K) $ have the same rank.
Let $$a:=| A_K(K)/\Pic^0(X_K)|.$$ 
We find in \cite{LLR1}, Proof of 4.6, based on the proofs of 2.3 and 2.5 in \cite{Gon}, that $a$ divides $(\prod {\delta_v'})/\lcm(\delta_v')$.
We are now ready to state the main result of this corrigendum.

\medskip
\noindent
{\bf Corrected Theorem 4.3.} 
{\it
Let $X/k$ and $f : X\to V$ be as above. 
Assume that $\Sha(A_K)$ and $\Br(X)$ are finite.
Then the equivalence of the Artin-Tate and Birch-Swinnerton-Dyer 
conjectures holds exactly when 
\begin{equation} \label{eq1}
|\Sha(A_K)| \prod_v \delta_v \delta'_v = 
a^2 \delta^2  |\Br(X)|.
\end{equation}
}
\medskip

The statement of Theorem 4.3 of \cite{LLR1} unfortunately omits the factor $a^2$
in the above formula. This omission leads to the following change in Corollary 4.7 of \cite{LLR1}. The paragraph after Corollary 4.7 in \cite{LLR1} can now be completely omitted. 

\medskip
\noindent
{\bf Corrected Corollary 4.7} {\it Assume that $\Sha(A_K)$ and  ${\rm Br}(X)$  are finite. Then
the conjectures of Artin--Tate and Birch--Swinnerton-Dyer are equivalent if and only if $\delta a = \delta' bc \epsilon$. }
\medskip

Theorem 4.3 in \cite{LLR1} is used in the proof of Corollary 3 of \cite{LLR2}.
The corrected version of Theorem 4.3 can be used in that proof to produce exactly the same result.
We restate below Corollary 3 of \cite{LLR2} with the correct formula relating the orders of $\Sha(A_K)$ and $\Br(X)$.

\medskip
\noindent
{\bf Corrected Corollary 3.} {\it Let $f : X \to V$ be as above. Assume that for some prime $\ell$,
the $\ell$-part of the group ${\rm Br}(X)$ or of the group $\Sha(A_K)$ is finite. Then
$|\Sha(A_K)| \prod_v \delta_v \delta'_v = a^2  \delta^2  |\Br(X)|$, and $|{\rm Br}(X)|$ is a square.}

\medskip

\end{section}
\begin{section}{Proof of the Corrected Theorem 4.3}

We follow closely the paper \cite{Gor} of Gordon, and indicate below every change that needs to be made to the statements in \cite{Gor} 
to obtain a complete proof of Formula \eqref{eq1}. 

\begin{emp} It may be of interest to first quickly indicate why the change in the formula occurs as a `square'. This fact is essential for the proof of Corollary 3 in \cite{LLR2} to remain correct.
The conjectures of Birch--Swinnerton-Dyer and of Artin--Tate require the explicit computation on one hand 
of the determinant of the height pairing on the lattice $A_K(K)/A_K(K)_{\rm tors}$, and on the other hand of
the determinant of the intersection pairing on the free part $NS(X)/NS(X)_{\rm tors}$ of the N\'eron--Severi group $NS(X)$. For this, it suffices to construct 
explicit bases for sublattices of finite index in these lattices (see, e.g., \ref{mainsublattice}, \ref{emp.pairing}), and the following well-known lemma then introduces `squares' in the formuli.
\end{emp}

\begin{lemma} \label{discr}
Let $\Lambda$ be a free abelian group of finite rank $n$, and let $\Lambda' \subseteq \Lambda$ be a sublattice of finite index $[\Lambda: \Lambda']$. 
Let $B: \Lambda \times \Lambda \to {\mathbb R}$ be a bilinear form. Consider a basis $\lambda_1,\dots,\lambda_n$ for $\Lambda$, and a basis
$\lambda'_1,\dots,\lambda'_n$ for $\Lambda'$. Let $d:= \det((B(\lambda_i,\lambda_j))_{1 \leq i,j\leq n})$, and similarly, let $d':= \det((B(\lambda'_i,\lambda'_j))_{1 \leq i,j\leq n})$. Then 
$$d'=[\Lambda: \Lambda']^2 d.$$
\end{lemma}

\begin{emp} \label{dvert} 
We introduce below a finite group $E$. This group
is claimed in \cite{Gor}, Lemma 4.3, to be always trivial, 
but the proof provided in \cite{Gor} is unfortunately incorrect (in the last paragraph, the computation of
$\pi^*C$ is wrong). 
This group will appear in two quotients of the filtration of $NS(X)$ introduced in \ref{Go4.5}. 
The final index discussed in \ref{Go4.6} however does not depend on $|E|$.

We follow below the notation in \cite{Gor} on page 177.   Let $\overline{k}$ denote an algebraic closure of $k$, and for any $k$-scheme $S$,  set 
as usual $\overline{S}:=S \times_k \overline{k}$. The natural map $\overline{X} \to X$ defines an injection $\Div(X) \to \Div(\overline{X})$ which is compatible with the intersection pairings $( \ , \ )_X$ and $( \ , \ )_{\overline{X}}$.
We identify $\Div(X)$ with its image in $\Div(\overline{X})$. Similarly, we use the maps $f: X \to V$ and $\overline{f}: \overline{X} \to \overline{V}$ to identify $\Div(V)$ and  $\Div(\overline{V})$ with their images in $\Div(X)$ and  $\Div(\overline{X})$, respectively.
Let us now define some natural subgroups of $\Div(\overline{X})$.

First, $\Div_{\rm vert}(\overline{X})$ is the subgroup generated
by the irreducible curves $C$ on $\overline{X}$ for which
$\overline{f}(C)$ is a single point. We denote by ${\rm
  Div}^0(\overline{X})$ the subgroup generated by the irreducible
curves $C$ on $\overline{X}$ which are algebraically equivalent to zero. Finally, let $\Div^0(\overline{V})$ 
denote the image in  $\Div(\overline{X})$ of 
 the subgroup of divisors on $\overline{V}$ algebraically  equivalent to zero. The subgroup  $\Div^0(\overline{V})$
 is the set of all divisors of the form $\sum_v a_v X_v$, where $X_v$ is the fiber over $v \in \overline{V}$ and $\sum_v a_v=0$.
The intersection of $\Div(X)$ with the subgroup $\Div^0(\overline{V})$, resp.\ with $\Div^0(\overline{X})$ or $\Div_{\rm vert}(\overline{X})$, 
is denoted by $\Div^0(V)$, resp.\ by $\Div^0(X)$, or  $\Div_{\rm vert}(X)$.

It is clear that $\Div^0(V)$  is contained in $\Div^0(X) \cap \Div_{\rm vert}(X)$. We let 
$$E:=  \frac{\Div^0(X) \cap \Div_{\rm vert}(X)}{\Div^0(V)}.$$

For $v \in V$, write $X_v=\sum_a p_{va} X_{va}$ with $X_{va}/k(v)$ irreducible of multiplicity $p_{va}$, and set $d_v:=\gcd_v(p_{va})$. 
The integer $d_v$ is called the multiplicity of the fiber $X_v$, and when $d_v>1$, $X_v$ is called a {\it multiple} fiber. 
Clearly $\frac{1}{d_v}X_v \in \Div(X)$.

If $W \in \Div^0(X) \cap \Div_{\rm vert}(X)$, then $W$ is numerically equivalent to zero, and so $(W \cdot X_{va})_X=0$ for all $X_{va}$. 
It follows from the fact that $\frac{1}{d_v}X_v$ generates the kernel of the intersection matrix associated with the fiber $X_v$
that $W = \sum_v m_v(\frac{1}{d_v}X_v)$ for some integers $m_v$.  Since $(W \cdot \Omega)_X=0$ for any horizontal divisor $\Omega$ on $X$, we find that
 $\sum_v (m_v/d_v) \deg_k v=0$. 
Hence for any $W\in \Div^0(X)\cap \Div_{\rm vert}(X)$, we have 
$W\in \Div^0(V)$ if and only if 
$m_v\in d_v\Z$ for all $v$. This implies that $E$ is isomorphic to a subgroup of 
$\oplus_v \Z/d_v\Z$. Let $\Delta:=\lcm_v(d_v)$. Then $E$ is killed by
$\Delta$ and $|E|$ divides $\prod d_v$.
\end{emp}

Let now $D_\ell(\overline{X})$ denote the subgroup of divisors in $\Div(\overline{X})$ that are linearly equivalent to zero. 
Set $D_\ell(X):= D_\ell(\overline{X}) \cap \Div(X)$. Let $\Pic^0_{X/k}$ and 
$\Pic^0_{V/k}$ denote the Picard schemes of $X/k$ and $V/k$, respectively.
($\Pic^0_{V/k}$ is nothing but the Jacobian of $V/k$.) The scheme 
$\Pic^0_{X/k}$ might not be reduced, and we 
denote by $\Pic^0_{X/k, {\rm red}}$ the (reduced) abelian variety associated with $\Pic^0_{X/k}$. 
We have
$$\Pic^0_{X/k, {\rm red}}(k) = \Div^0(X)/D_\ell(X) 
\text{ and } \Pic^0_{V/k}(k) = \Div^0(V)/D_\ell(V)$$
because $\Br(k)$ is trivial.

\begin{lemma} Keep the above notation. Then 
\begin{enumerate}[\rm a)]
\item We have $(\Div^0(X) \cap \Div_{\rm vert}(X)) \cap (\Div^0(V)+D_{\ell}(X))= \Div^0(V)$.
\item We have a natural injection 
$$E \longrightarrow \Pic^0_{X/k, {\rm red}}(k)/\Pic^0_{V/k}(k)$$
given explicitly as 
\begin{multline*} 
\frac{\Div^0(X) \cap \Div_{\rm vert}(X)}{\Div^0(V)}=  
\frac{ 
\Div^0(X) \cap (\Div_{\rm vert}(X)+D_{\ell}(X)) }{ 
\Div^0(V)+D_{\ell}(X)
 } 
\longrightarrow  \\ 
\longrightarrow \frac{\Div^0(X)}{\Div^0(V)+D_{\ell}(X)}. \hskip 30pt 
\end{multline*}
\end{enumerate}
\end{lemma}
\proof The proof of b) follows immediately from a). To prove Part a), it suffices to prove that 
$$ \Div_{\rm vert}(X) \cap (\Div^0(V)+D_{\ell}(X))= \Div^0(V).$$
If $D \in \Div_{\rm vert}(X) \cap (\Div^0(V)+D_{\ell}(X))$, then $D \in \Div_{\rm vert}(X) \cap \Div^0(X)$.
 As noted in \ref{dvert}, we can then write $D = \sum_v r_vX_v$ for some rational numbers $r_v$ with $\sum_v r_v \deg(v)=0$. 
On the other hand, by hypothesis, $D={\rm div}(f) + D_0$ for some $f \in k(X)^*$ and $D_0 \in  \Div^0(V)$. 
Since $k$ is finite, some multiple of $D_0$ is linearly equivalent to
zero. Thus, for some positive integer $m$, $mD={\rm div}(f^mh)$ for
some $h   \in k(V)^*$. 
Since $mD = \sum_v mr_vX_v\in \Div^0(V)$, we find that some positive multiple $n$ of $mD$ is of the form ${\rm div}(h')$ for some $h' \in  k(V)^*$.
Hence, $f^{mn} \in k(V)^*$. Since we assume that the generic fiber of
$X \to V$ is 
geometrically integral, we find that $f \in k(V)^*$. Thus $D \in \Div^0(V)$.
\qed

\medskip
We stray here a little bit from the notation used by \cite{Gor}, and we define $B/k$ to be the quotient
abelian variety $B:= \Pic^0_{X/k, {\rm red}} / \Pic^0_{V/k}$.
Since $k$ is finite, we have 
$$B(k):= \Pic^0_{X/k, {\rm red}}(k) / \Pic^0_{V/k}(k).$$
For use in the proof of \ref{Go4.5} (iv), let us note that
\begin{equation} \label{eqB/E}
\frac{B(k)}{E} = \frac{\Div^0(X)}{(\Div^0(X) \cap \Div_{\rm vert}(X))+D_{\ell}(X)}.
\end{equation}

\begin{remark} In \cite{Gor}, just before Proposition 4.4 on page 180, $B/k$ is defined to be the $K/k$-trace of $A_K/K$. 
Then Proposition 4.4 asserts that the $K/k$-trace of $A_K/K$ is an abelian variety which is purely inseparably isogenous to the quotient abelian variety 
$ \Pic^0_{X/k, {\rm red}} / \Pic^0_{V/k}$. The proof of Proposition 4.4 in \cite{Gor} uses the fact that $a=1$.  We refer the reader to \cite{Con} for the definition and existence of the $K/k$-trace of $A_K/K$. When $k$ is algebraically closed, we find in \cite{Shi}, Theorem 2, a theorem of Raynaud
which asserts that the $K/k$-trace of $A_K/K$ is $k$-isomorphic to $ \Pic^0_{X/k, {\rm red}} / \Pic^0_{V/k}$ when  $f: X \to V$ does not have any multiple fibers (i.e., $d_v=1$ for all $v$). The notion of $K/k$-trace is not needed in this corrigendum, and we do not use Proposition 4.4 in \cite{Gor}.
\end{remark}

Let $\Div_0(\overline{X})$ denote the subgroup of $\Div(\overline{X})$ generated by the irreducible curves which intersect each complete vertical fiber $X_v$ with total intersection multiplicity zero. We let $\Div_0(X):= \Div_0(\overline{X})\cap \Div(X)$.
 Let $\Omega  \in \Div(X)$ be a horizontal divisor of degree $\delta$, where $\delta$ is the index of $X_K$ over $K$.
 In the following modified version of Lemma 4.2 in \cite{Gor}, 
 the group $A_K(K)$ has now been replaced by $\Pic^0(X_K)$.

 \begin{lemma} \label{Lemma4.2} {\rm (see Lemma 4.2 in \cite{Gor})}
 There  are natural isomorphisms of groups
 $$\frac{\Div(X)}{(\Div_{\rm vert}(X) \oplus {\mathbb Z}\Omega)+D_\ell(X)} \longrightarrow 
 \frac{\Div_0(X)}{\Div_{\rm vert}(X) +D_\ell(X)}  \longrightarrow \Pic^0(X_K).$$
 \end{lemma}
 \proof Same as in \cite{Gor}, replacing when necessary $A_K(K)$  by $\Pic^0(X_K)$. \qed

\begin{emp} \label{mainsublattice} Let $NS(X):= \Div(X)/\Div^0(X)$. Let us now introduce further notation needed to  define below
the completely explicit subgroup ${\mathcal N}_0$ of $NS(X)$.
\begin{enumerate}[\rm (a)]
\item
Let $r$ be the rank of $A_K(K)$, 
 and let $\{ \alpha_1, \dots, \alpha_r\}$ be a basis of the lattice $\Pic^0(X_K)/ \Pic^0(X_K)_{\rm tors}$.
 Choose divisors ${\mathcal A}_1$, $\dots$, ${\mathcal A}_r$ in $\Div(X)$ such that for each $i$, the class in $\Pic^0(X_K)$
 of the restriction of ${\mathcal A_i}$ to the generic fiber $X_K$ is $\alpha_i$. 
 For the later purpose of computing the global height pairing $\langle \alpha_i, \alpha_j\rangle$ as in \ref{globalheight},
 we assume also that we have chosen the divisors ${\mathcal A}_1$, $\dots$,  ${\mathcal A}_r$, such that the supports of the restrictions of ${\mathcal A_i}$ and ${\mathcal A_j}$ to the generic fiber $X_K$ are pairwise disjoint when $i \neq j$.
\item Since $X_K/K$ has index $\delta$, choose a divisor $\sum_i s_ix_i$ in $\Div(X_K)$ such that $\sum_i s_i \deg_K(x_i)=\delta$.
Let $\overline{x_i}$ denote the closure of $x_i$ in $X$, and set $\Omega:=\sum_i s_i \overline{x_i}$ in $\Div(X)$.
\item Since $V/k$ is geometrically integral, its index $\delta(V/k)$ is equal to $1$. Choose a divisor $\sum_j t_j v_j$ in $\Div(V)$ such that 
$\sum_j t_j \deg_k(v_j)=1$. Let $F:= \sum_j t_j X_{v_j}$ in $\Div(X)$. This definition agrees with \cite{Gor}, 4.6, when $X_K$ has a $k$-rational point and the complete fiber in 4.6 is chosen to be above a $k$-rational point.
\item For each $v\in V$, write the fiber $X_v$ as $X_v=\sum_{a=1}^{h(v)} p_{va} X_{va}$, where the components $X_{va}$ are irreducible.
For each closed point $v \in V$ such that $X_v$ is reducible, 
consider the set $\{ X_{va}, a>1, v \in V\}$ of irreducible divisors in $\Div(X)$. 
\end{enumerate}

We let ${\mathcal N}_0$ denote the subgroup of $NS(X)$ generated by 
$NS(X)_{\rm tors}$ and the classes of $\{ {\mathcal A}_1, \dots,
{\mathcal A}_r\}$, $\Omega$, $F$, and $\{ X_{va}, a>1, v \in V\}$. We
will compute the index of ${\mathcal N}_0$ in $NS(X)$ in Proposition~\ref{Go4.6}. 
\end{emp}

Denote by 
$S_1$ the set of closed points $v \in V$ such that  $X_v$ is
reducible. Let $S_2$ denote the set of closed points $v \in V$ such that $X_v$ is irreducible but not reduced. Set $\Sigma:=S_1 \sqcup S_2$. Let $S_3$ denote the set of $v \in V$ such that $X_v$ is integral but not geometrically integral. 

The set $\Sigma$ is finite, and thus we have  
\begin{equation} \label{eq:defQ}
Q:=\frac{\Div_{\rm vert}(X)}{\Div(V)}= \frac{\oplus_v (\oplus_a {\mathbb Z}X_{va})}{\oplus_v {\mathbb Z}X_v}
= \oplus_{v\in \Sigma} \left( \frac{\oplus_a {\mathbb
      Z}X_{va}}{{\mathbb Z}X_{v}}\right).
\end{equation} 
Define $NS(X)_{\rm vert}$ to be the image in $NS(X)$ of the subgroup
$\Div_{\rm vert}(X)$ of $\Div(X)$. Let $[\Omega]$ denote the class of $\Omega $ in $NS(X)$. 
It is clear that $NS(X)_{\rm vert} \cap {\mathbb Z}[\Omega] = (0)$,
and we write 
$${\mathcal N}:=NS(X)_{\rm vert} \oplus {\mathbb Z}[\Omega].$$
We may now state a modified version of Proposition 4.5 in \cite{Gor}, 
where the group $E$ occurs in two different factors.

\begin{proposition} \label{Go4.5} {\rm (see Proposition 4.5 in \cite{Gor})}
The group $NS(X)$ has a filtration by subgroups
$$0 \subseteq f^*NS(V) \subseteq NS(X)_{\rm vert} \subseteq {\mathcal N} \subseteq NS(X)$$
with respective quotients ${\mathbb Z}$, $Q/E$, ${\mathbb Z}$, and 
$\Pic^0(X_K)/(B(k)/E)$.
\end{proposition}
\proof
(i) {\it The map $f^*: NS(V)\to NS(X)$ is injective}, and since $NS(V)$ is free of rank $1$, so is $f^*NS(V)$.

(ii) Let us first note 
that the natural map
$$E=\frac{\Div^0(X) \cap \Div_{\rm vert}(X)}{\Div^0(V)}
\longrightarrow Q=\frac{\Div_{\rm vert}(X)}{\Div(V)}$$
is injective because  
\begin{equation} \label{Eeq}
\Div^0(X) \cap \Div_{\rm vert}(X) \cap \Div(V) = \Div^0(V).
\end{equation}
Recall that $$NS(X)_{\rm vert}= \frac{\Div_{\rm vert}(X)}{\Div_{\rm vert}(X)\cap \Div^0(X)},$$
and consider the natural map 
$f^*\Div(V) \longrightarrow NS(X)_{\rm vert}$. 
This map has kernel $f^*\Div^0(V)$, by \eqref{Eeq}. Hence, we have an exact sequence
$$ 0 \to f^* NS(V) \longrightarrow  NS(X)_{\rm vert} \longrightarrow Q/E \longrightarrow 0.$$

(iii) By construction ${\mathcal N}/NS(X)_{\rm vert}=\Z [\Omega]\simeq
\Z$.

(iv)
As in Part (4) of the proof in \cite{Gor}, we have an exact sequence
\begin{multline*}
0 \longrightarrow \frac{\Div^0(X)}{\Div^0(X) \cap (\Div_{\rm vert}(X))
  + D_{\ell}(X))} \longrightarrow \\ 
\longrightarrow\frac{ \Div_0(X)}{\Div_{\rm vert}(X)
  + D_{\ell}(X)} \longrightarrow \frac{NS(X)}{{\mathcal N}}
\longrightarrow 0.
\end{multline*} 
The first term in this sequence is identified with $B(k)/E$ in \eqref{eqB/E} since $D_{\ell}(X) \subseteq \Div^0(X)$. 
The middle term is identified with $\Pic^0(X_K)$ in \ref{Lemma4.2}. 
 We thus have an isomorphism 
$$NS(X)/{\mathcal N} \longrightarrow \frac{\Pic^0(X_K)}{B(k)/E}.$$
\qed

\begin{proposition} \label{Go4.6} {\rm (see Proposition 4.6 in \cite{Gor})}
Let $\N_0\subseteq NS(X)$ be as in {\rm \ref{mainsublattice}}.
Then the quotient 
$NS(X)/\N_0$ is finite with 
$$|NS(X)/\N_0|=\frac{ |\Pic^0(X_K)_{\rm tors}| }{|B(k)|} \cdot \frac{\prod_{v \in \Sigma} p_{v1}}{|NS(X)_{\rm tors}|}.
$$
\end{proposition}

\proof 
Let ${\mathcal N}'$ be the subgroup of $NS(X)$ generated by the
classes of $\Omega, F,$ and $ X_{va}$ for $a>1$ and $h(v)>1$, so that ${\mathcal N}' \subseteq {\mathcal N}_0 $. 
Recall that $\N:=NS(X)_{\rm vert} \oplus {\mathbb Z}[\Omega]$, so that ${\mathcal N}' \subseteq \N $. 
We have two exact sequences 
\[
\xymatrix{& & 0\ar[d] & & \\
& & \N/\N' \ar[d] & & \\ 
0\ar[r] & \N_0/\N' \ar[r] & A':=NS(X)/\N' \ar[r]\ar[d] & NS(X)/\N_0 \ar[r] & 0 \\ 
& & P:=NS(X)/\N \ar[d] & & \\ 
& & 0 & &} 
\] 

Let us start by computing the order of ${\mathcal N}/{\mathcal N}'$. Write ${\mathcal N}''$
for the subgroup of ${\mathcal N}'$ generated by the classes of $F$, and $ X_{va}$, $a>1$ for all $v$ with $h(v)>1$.
Then ${\mathcal N}'' \subseteq NS(X)_{\rm vert}$ and 
$\N'=\N'' \oplus {\mathbb Z}[\Omega]$. 
It follows that 
$$
\frac{{\mathcal N}}{{\mathcal N}'} = \frac{NS(X)_{\rm vert}}{{\mathcal N}''} = \frac{NS(X)_{\rm vert}/f^*NS(V)}{({\mathcal N}''+f^*NS(V))/f^*NS(V)}.$$
The numerator of the group on the right is identified with $Q/E$ in \ref{Go4.5}. One  checks that ${\mathcal N}''\cap f^*NS(V)= {\mathbb Z}[F]$. 
With the group $Q$ identified as in \eqref{eq:defQ}, let $Q'$ denote the subgroup of $Q$ generated by the classes of the
components $X_{va}$ with $a>1$ for all $v $ with
$h(v)>1$. Then the denominator in the above expression is equal to
$Q'$ and it is clear that $Q/Q'$ is isomorphic to $\prod_{v \in \Sigma} {\mathbb Z}/p_{v1}{\mathbb Z}$. Since $Q'$ is torsion free and $E$ is torsion, 
we find that 
$${\mathcal N}/{\mathcal N}'\simeq (Q/E)/Q'\simeq Q/(Q'+E),$$ 
so that ${\mathcal N}/{\mathcal N}'$ is finite, of order $(\prod_{v \in \Sigma}p_{v1})/|E|$. 

Recall now from \ref{Go4.5} that $P\simeq \Pic^0(X_K)/(B(k)/E)$. 
Since $B(k)/E$ is finite, we find that 
\begin{equation} \label{P0}
|P_{\rm tors}|= |\Pic^0(X_K)_{\rm tors}|/|B(k)/E|,
\end{equation}
and we also have a canonical isomorphism
\begin{equation} \label{P}
\Pic^0(X_K)/  \Pic^0(X_K)_{\rm tors}  \longrightarrow P/P_{\rm tors}.
\end{equation}
Since the group ${\mathcal N}/{\mathcal N}'$ is finite, we find that
\begin{equation} \label{A0}
|A'_{\rm tors}| = |{\mathcal N}/{\mathcal N}'| \cdot |P_{\rm tors}|
\end{equation}
and that 
\begin{equation} \label{A}
A'/A'_{\rm tors} \longrightarrow P/P_{\rm tors}
\end{equation}
is an isomorphism.

By construction, the classes of the restrictions of ${\mathcal A}_1,\dots,{\mathcal A}_r$ to the generic fiber are a basis of $\Pic^0(X_K)/  \Pic^0(X_K)_{\rm tors}$. 
Using the isomorphisms \eqref{P} and \eqref{A}, we find that the classes of ${\mathcal A}_1,\dots,{\mathcal A}_r$ are a basis of $A'/A'_{\rm tors}$.
This implies that $NS(X)/\N_0$ is torsion and that 
$$ 0 \longrightarrow (\N_0/\N')_{\rm tors} \longrightarrow A'_{\rm tors} \longrightarrow NS(X)/\N_0 \longrightarrow 0$$
is exact. It is clear that 
$$\N_0=(\langle [{\mathcal A}_1], \dots, [{\mathcal
    A}_r]\rangle + NS(X)_{\rm tors}) \oplus \N'.$$ 
It follows that 
$$NS(X)/\N_0= \frac{A'_{\rm tors}}{NS(X)_{\rm tors}}.
$$
The desired formula for the index follows from \eqref{P0} and \eqref{A0}.
\qed

\begin{emp} \label{emp.pairing} Let $\overline{\N}_0$ be the image of $\N_0$ in the
lattice   $NS(X)/NS(X)_{\rm tors}$. 
The computation of the discriminant of the intersection pairing on the
sublattice 
$\overline{\N}_0$ is done exactly as in 
Proposition 5.1 of \cite{Gor}, and the formula obtained is the same. The only difference now is that 
the discriminant of the height pairing $|\det\left< \alpha_i, \alpha_j\right>|$ that appears   in the formula is the discriminant
for the height pairing on $\Pic^0(X_K)/ \Pic^0(X_K)_{\rm tors}$, and not anymore on $A_K(K)/ A_K(K)_{\rm tors}$.
Let $a_f$ denote the index of $\Pic^0(X_K)/ \Pic^0(X_K)_{\rm tors}$ in $A_K(K)/ A_K(K)_{\rm tors}$. As indicated in Lemma \ref{discr}, the two discriminants
differ by a factor $a_f^2$.

Similarly, the discriminant of the intersection pairing on $\overline{\N}_0$
differs from the discriminant of the intersection pairing 
on the full lattice $NS(X)/NS(X)_{\rm tors}$ by the square of the index
$$ \frac{ |\Pic^0(X_K)_{\rm tors}| }{|B(k)|} \cdot \frac{\prod_{v \in \Sigma} p_{v1}}{|NS(X)_{\rm tors}|}
$$
obtained in \ref{Go4.6}. This index is exactly the same as the one obtained \cite{Gor}, except that in \cite{Gor}, 
the term $|\Pic^0(X_K)_{\rm tors}|$ is replaced by $|A_K(K)_{\rm tors}|$. Let $a_{\rm tors}:=|A_K(K)_{\rm tors}/\Pic^0(X_K)_{\rm tors}|$. 
We have $a=a_f a_{\rm tors}$, and we find that the final discrepancy is a factor of $a^2$.

\end{emp}

\begin{remark} \label{globalheight}
We supply in this remark some references for an important result stated just before Proposition 5.1 of \cite{Gor}, and needed in its proof.
Let $\alpha, \beta$ in $\Pic^0(X_K)/ \Pic^0(X_K)_{\rm tors}$.
The global height pairing 
$\left< \alpha, \beta \right>$ can be computed as a sum of local contributions $\sum_v \left< \alpha, \beta \right>_v$ (see, e.g., \cite{Gro}, (4.6)).
Each local contribution  
can be expressed as a local intersection number $\left< \alpha, \beta \right>_v= - (\alpha , \beta)_v \log(|k(v)|)$ (see, e.g., \cite{Gro}, (3.7)),
where the contribution $(\alpha , \beta)_v$ is the value of N\'eron's pairing at $v$ on $\alpha$ and $\beta$.
Let $A, B \in \Div(X) \otimes {\mathbb Q}$ be two divisors whose restrictions to $X_K$ are in $\Div(X_K)$ and equal the classes $\alpha$ and $\beta$, respectively, and have disjoint supports. Assume in addition that $(A \cdot X_{va})_X=0$ for all $v$ and all $a$. 
Then $(\alpha , \beta)_v = (A \cdot B)_v$, where $(A \cdot B)_v$ denotes the contribution of the points in $X_v$ in the intersection number $(A \cdot B)_X$ 
(see, e.g., \cite{BLo}, 4.3, or \cite{Pep}, 2.2).
One then obtains that $\left< \alpha, \beta \right>= - (A \cdot B)_X \log(|k|)$.
\end{remark}

\begin{emp}	 
We recall below the formula of Gordon
found in the middle of page 196 in \cite{Gor}. This formula is claimed to hold exactly when the Birch--Swinnerton-Dyer conjecture is equivalent to the Artin--Tate conjecture. This  claim is incorrect when $a>1$. In \cite{Gor}, page 169, the integer  $\alpha $ appearing below is defined to be the index $\delta$. 

\begin{equation} \label{eqchi2}
|\Sha(A_K)| \prod_v d_v^2 \epsilon_v = 
 \alpha^2  |\Br(X)|.
\end{equation}

This formula in \cite{Gor} is misleading, as the term $\epsilon_v$ is only introduced in the statement of Proposition 5.5 of \cite{Gor} when $v \in S_1$,
but the formula (6.2) in \cite{Gor}, from which  \eqref{eqchi2} above is derived, involves a product over a set $S$ 
(defined on page 165 of \cite{Gor})
which contains $S_1$, but which might also contain $S_2$ and $S_3$ (notation introduced in \ref{mainsublattice}). Let us therefore state below the correct formula \eqref{eqchi2correct} that can be inferred from Gordon's work and which should be substituted for 
 \eqref{eqchi2}. 
 
 Let ${\mathcal A}_v/{\mathcal O}_{K_v}$ denote the N\'eron model of $A_{K_v}/K_v$. 
 Let $\Phi_v/k(v)$ denote the group of components of the special fiber of ${\mathcal A}_v$. When $v \in S_2 \sqcup S_3$, the fiber $X_v$ is irreducible, 
 say $X_v= d_v \Gamma_v$ for some irreducible curve $\Gamma_v/k(v)$. Let $q_v$ denote the degree over $k(v)$ of the algebraic closure of $k(v)$ in the function field of $\Gamma_v/k(v)$. 
 It follows from the fact that $k(v)$ is a finite field that $\delta_v = d_v q_v$. Note that if $v \notin S_1 \sqcup S_2 \sqcup S_3$, then $\delta_v=\delta'_v =1$.
 Then Gordon's arguments, along with the removal of the hypothesis that $X\to V$ be cohomologically flat in dimension $0$ in \cite{LLR1} 
 and the corrections given in this corrigendum, give the following formula.

 \begin{equation} \label{eqchi2correct}
|\Sha(A_K)| \left( \prod_{v \in S_1} d_v^2 \epsilon_v\right)\left(  \prod_{v \in S_2 \sqcup S_3} d_v^2 |\Phi_v(k(v))| q_v\right) =  a^2 \delta^2  |\Br(X)|.
\end{equation}
The formula can be turned into Formula \eqref{eq1} as we did in the proof of Theorem 4.3 in \cite{LLR1}, using Theorem 1.17 of \cite{B-L}. 
For instance, when $v \in S_2 \sqcup S_3$, this theorem shows that $|\Phi_v(k(v))|= \delta_v'/d_v$. Since it follows from the adjunction formula
that $d_vq_v$ divides $g-1$ in this case,  Theorem 7 in \cite{Lic} shows that $\delta_v=\delta_v'$. It follows that  $d_v^2 |\Phi_v(k(v))| q_v= \delta_v \delta_v'$, as desired, and Formula \eqref{eq1} is established.
\end{emp}

\end{section}

\end{document}